# INVERSE PROBLEMS FOR REGULAR VARIATION OF LINEAR FILTERS, A CANCELLATION PROPERTY FOR σ-FINITE MEASURES AND IDENTIFICATION OF STABLE LAWS

By Martin Jacobsen,[1] Thomas Mikosch,[1] Jan Rosiński[2]
and Gennady Samorodnitsky[3]

*University of Copenhagen, University of Tennessee and Cornell University*

In this paper, we consider certain σ-finite measures which can be interpreted as the output of a linear filter. We assume that these measures have regularly varying tails and study whether the input to the linear filter must have regularly varying tails as well. This turns out to be related to the presence of a particular cancellation property in σ-finite measures, which in turn, is related to the uniqueness of the solution of certain functional equations. The techniques we develop are applied to weighted sums of i.i.d. random variables, to products of independent random variables, and to stochastic integrals with respect to Lévy motions.

**1. Introduction.** Regular variation is one of the basic concepts which appear in a natural way in different contexts of applied and theoretical probability. Its use in the universe of applied probability models is well spread, as illustrated in the encyclopedic treatment in [4]. We will discuss only regular variation in one dimension though some of the questions addressed below have natural counterparts for multivariate regular variation. Recall that a random variable $Z$ is said to have a *regularly varying (right) tail with exponent* $\alpha > 0$ if its law $\mu$ satisfies the relation

$$(1.1) \qquad \mu(x, \infty)(= P(Z > x)) = x^{-\alpha} L(x) \qquad \text{for } x > 0,$$

Received December 2007; revised March 2008.

[1]Supported in part by the Danish Research Council (FNU) Grants 21-04-0400 and 272-06-0442.

[2]Supported in part by NSA Grant MSPF-07G-126.

[3]Supported in part by NSA Grant MSPF-05G-049, the ARO Grant W911NF-07-1-0078 at Cornell University and a Villum Kann Rasmussen grant at the University of Copenhagen.

where $L$ is a slowly varying (at infinity) function. Here and in what follows, we write for convenience $\nu(a, b] = \nu((a, b])$, $\nu[a, b] = \nu([a, b])$, and so forth, for any measure $\nu$ on $\mathbb{R}$ and $-\infty \leq a < b \leq \infty$. It is common to say that $Z$ is *regularly varying with index* $\alpha$. Finally, we use relation (1.1) as the defining property of a regularly varying right tail for any measure $\mu$ on $\mathbb{R}$ as long as $\mu(a, \infty) < \infty$ for some $a > 0$ large enough. If (1.1) holds with $\alpha = 0$, it is common to use the alternative terms *slowly varying* or *regularly varying with exponent* 0.

Regular variation has the important property that it is preserved under a number of linear operations and transformations which are often applied in probability theory. In other words, if the stochastic input to a linear filter has a regularly varying tail, then (under mild assumptions), the output from the filter also has a regularly varying tail. In this paper, we are interested in the *inverse problem*: suppose that the output from a linear filter is regularly varying with index $\alpha \geq 0$. When may we conclude that the input to the filter is regularly varying with the same index? The question is of obvious interest in a variety of filtering problems. To make it clear what we have in mind, we start with three examples.

EXAMPLE 1.1 (Weighted sums). Let $Z_1, Z_2, \ldots$ be i.i.d. random variables, and $\psi_1, \psi_2, \ldots$, nonnegative weights. If a generic element $Z$ of the sequence $(Z_j)$ is regularly varying with exponent $\alpha > 0$, then under appropriate conditions on the coefficients, the infinite series

$$(1.2) \qquad X = \sum_{j=1}^{\infty} \psi_j Z_j$$

converges with probability 1, and

$$(1.3) \qquad \lim_{x \to \infty} \frac{P(X > x)}{P(Z > x)} = \sum_{j=0}^{\infty} \psi_j^{\alpha}.$$

In general, the conditions on the coefficients $\psi_j$ ensuring summability of $X$ can be affected by the left tail of $Z$ as well. We refer to [21] for the most general conditions of this kind which are close to the necessary and sufficient conditions for summability of $X$. The equivalence (1.3) is always true for finite sums and regularly varying $Z$.

Weighted sums of type (1.2) arise naturally as the marginals of (perhaps infinite) moving average (or linear) processes. They constitute the backbone of classical linear time seriesanalysis; see, for example, [7]. Over the last 25 years, regular variation for these models has been used very effectively to a large extent through extreme value theory and point process techniques; see, for example, [8, 9, 10] for some pioneering work in time series analysis for heavy-tailed linear processes; compare [15], Chapter 7, for a survey of related results.

The research on infinite weighted sums as in (1.2) with slowly varying noise variables does not appear to be far advanced, but in the case of *finite sums* (i.e., if only finitely many of the coefficients $\psi_j$ are different from zero) relation (1.3) is still well known to hold; see, for example, [17].

EXAMPLE 1.2 (Products). Let $Z$ be a random variable that is regularly varying with exponent $\alpha \geq 0$, independent of another random variable $Y > 0$, and write $X = YZ$. If the tail of $Y$ is light enough, then the tail of $X$ is also regularly varying with exponent $\alpha$. If, for example, $EY^{\alpha+\varepsilon} < \infty$ for some $\varepsilon > 0$, then the tail equivalence

$$(1.4) \qquad \lim_{x\to\infty} \frac{P(X > x)}{P(Z > x)} = EY^{\alpha}$$

holds; see [6]. It is also well known that $X$ is regularly varying with index $\alpha \geq 0$ if both $Z$ and $Y$ are regularly varying with index $\alpha$; see [13], and for variations on this theme, [11].

EXAMPLE 1.3 (Stochastic integrals). Let $(M(s))_{s\in\mathbb{R}}$ be a Lévy process with Lévy measure $\eta$ (see [28] for the general theory of Lévy processes), and $g: \mathbb{R} \to \mathbb{R}_+$ a measurable function. Under certain integrability assumptions on $g$, the random variable

$$(1.5) \qquad X = \int_{\mathbb{R}} g(s) M(ds)$$

is well defined; see [23]. If the Lévy measure $\eta$ has a regularly varying tail with exponent $\alpha \geq 0$, then the integral $X$ itself is regularly varying and

$$(1.6) \qquad \lim_{x\to\infty} \frac{P(X > x)}{\eta(x,\infty)} = \int_{\mathbb{R}} [g(s)]^{\alpha}\, ds,$$

see [26]. In the case $\alpha = 0$, this last relation requires the function $g$ to be supported by a set of a finite Lebesgue measure.

Stochastic integrals of type (1.5) naturally appear as the marginals of continuous time moving average processes, one of the major classes of stochastic models and other infinitely divisible processes.

It is possible to view the three examples above as demonstrating the preservation of regular variation by a linear filter. In fact, the random variable $X$ in Example 1.1 (or, even more clearly, an entire linear process) is the output of the linear filter defined by the coefficients $\psi_j$, $j = 1, 2, \ldots$ in (1.2). The product $X$ in Example 1.2 is a random linear transformation of the input $Z$. The stochastic integral (1.5) in Example 1.3 can be viewed as a linear transformation of the increments of the Lévy process $M$. These examples form the heart of this paper, and for those examples the inverse problem, in which we are presently interested, can be formulated very explicitly. In Example 1.1, suppose that the series (1.2) converges and the sum $X$ is regularly varying. Does it necessarily follow that a typical noise variable $Z$ is also regularly varying? In Example 1.2, if the product $X = YZ$ is regularly varying and $Y$ is "appropriately small," is $Z$ regularly varying? In Example 1.3, if the stochastic integral (1.5) is regular varying, are so the increments of the Lévy process $M$, hence the tail of the Lévy measure $\eta$ (see [31])?

It turns out that the above inverse problems are related. Surprisingly, in the case $\alpha > 0$, they are all connected to two other, apparently unrelated, problems: that of functional equations of integrated Cauchy and Choquet–Deny type (see [24]), and that of the cancellation property of certain $\sigma$-finite measures. We set up these two latter problems in Section 2. Negative answers to the inverse problems on regular variation for the linear filters discussed above turn out to correspond to the existence of nontrivial solutions to these related functional equations as well as to the lack of the cancellation property of measures. In Section 2, we will also see that these phenomena are related to the existence of real zeros of certain Fourier transforms.

We apply these results to weighted sums in Section 3, to products of independent random variables in Section 4 and to stochastic integrals in Section 5. Interestingly, it appears that, "in most cases," the answer to the inverse problems on regular variation for linear filters is positive. This is emphasized by a numerical study in Section 3, where the case of a weighted sum with 3 weights is considered.

We will also consider the case of slow variation, that is, when $\alpha = 0$. In this case, the inverse problems are almost trivial in the sense that slow variation at the output of the linear filter implies slow variation at its input under very general conditions on the filter. Surprisingly, in contrast to the case $\alpha > 0$, there does not seem to be any connection between the inverse problems and the cancellation property of related measures.

In the final Section 6, we touch on a closely related but more special problem. It is well known (see, e.g., [17, 27]) that infinite variance $\alpha$-stable distributions are regularly varying with exponent $\alpha \in (0, 2)$. But given the output of a linear filter has such an $\alpha$-stable distribution, is it also true that the input to the filter is $\alpha$-stable? We indicate how to use our results for these distribution identification problems and give partial answers under suitable conditions on the filter.

**2. The cancellation property of measures and related functional equations.** Let $\nu$ and $\rho$ be two $\sigma$-finite measures on $(0, \infty)$. We define a new measure on $(0, \infty)$ by

$$(2.1) \quad (\nu \circledast \rho)(B) = \int_0^\infty \nu(x^{-1} B) \rho(dx), \qquad B \text{ a Borel subset of } (0, \infty).$$

Since $\nu \circledast \rho$ is the usual convolution of measures on the multiplicative group $(0, \infty)$, we call it the *multiplicative convolution* of the measures $\nu$ and $\rho$.

We say that a $\sigma$-finite measure $\rho$ has the *cancellation property* with respect to a family $\mathcal{N}$ of $\sigma$-finite measures on $(0, \infty)$ if for any $\sigma$-finite measures $\nu, \overline{\nu}$ on $(0, \infty)$ with $\overline{\nu} \in \mathcal{N}$,

$$(2.2) \qquad \nu \circledast \rho = \overline{\nu} \circledast \rho \quad \Rightarrow \quad \nu = \overline{\nu}.$$

The cancellation property for Lévy measures was considered in [1]. It was shown that certain Lévy measures $\rho$ have the cancellation property with respect to $\mathcal{N} = \mathcal{N}_\rho = \{\overline{\nu} : \overline{\nu} \circledast \rho \text{ is a Lévy measure}\}$ but a general characterization of such $\rho$'s remains unknown.

If $\mathcal{N} = \{\delta_1\}$, then (2.2) is known as the Choquet–Deny equation in the multiplicative form. The class of measures having the cancellation property with respect to $\{\delta_1\}$ can be determined by the well-studied Choquet–Deny theory, compare [24].

In this paper, we are only interested in the case when $\mathcal{N}$ consists of a single measure $\nu_\alpha$ with power density function, that is, $\nu_\alpha$ is a $\sigma$-finite measure on $(0, \infty)$ with density

$$\frac{\nu_\alpha(dx)}{dx} = \begin{cases} |\alpha| x^{-(\alpha+1)}, & \alpha \neq 0, \\ x^{-1}, & \alpha = 0. \end{cases}$$

If $\alpha \in (0, 2)$, $\nu_\alpha$ is the Lévy measure of an $\alpha$-stable law, whereas $\nu_0$ is the Haar measure for the multiplicative group $(0, \infty)$, and $\nu_{-1}$ corresponds to the Lebesgue measure. In the applications to the inverse problems for regular variation discussed in Section 1, only positive values of $\alpha$ are of interest; nonetheless, there are situations that require the cancellation property for nonpositive $\alpha$. An example of such situations is in Proposition 6.2 in Section 6.

THEOREM 2.1. *We assume that $\alpha \in \mathbb{R}$ and $\rho$ is a nonzero $\sigma$-finite measure such that*

$$(2.3) \qquad \int_0^\infty y^{\alpha - \delta} \vee y^{\alpha + \delta} \rho(dy) < \infty \qquad \text{for some } \delta > 0.$$

*Then the measure $\rho$ has the cancellation property with respect to $\mathcal{N} = \{\nu_\alpha\}$ if and only if*

$$(2.4) \qquad \int_0^\infty y^{\alpha + i\theta} \rho(dy) \neq 0 \qquad \text{for all } \theta \in \mathbb{R}.$$

*Furthermore, if the left-hand side of (2.4) vanishes for some $\theta_0 \in \mathbb{R}$, then for any real $a, b$ with $0 < a^2 + b^2 \leq 1$, the $\sigma$-finite measure*

$$(2.5) \qquad \nu(dx) := g(x) \nu_\alpha(dx)$$

*with*

$$g(x) := 1 + a\cos(\theta_0 \log x) + b\sin(\theta_0 \log x), \qquad x > 0, \tag{2.6}$$

*satisfies the equation* $\nu \circledast \rho = \nu_\alpha \circledast \rho$.

PROOF. We first consider the case $\alpha > 0$. Suppose that (2.4) holds. Note that the equation $\nu \circledast \rho = \nu_\alpha \circledast \rho$ can be rewritten as

$$\nu \circledast \rho = \|\rho\|_\alpha \nu_\alpha, \tag{2.7}$$

where by condition (2.3),

$$\|\rho\|_\alpha := \int_0^\infty y^\alpha \rho(dy) \in (0, \infty)$$

is the $\alpha$th moment of $\rho$. Let $\varepsilon > 0$ be such that $\rho(\varepsilon, \infty) > 0$. Then it follows from (2.7) that

$$\nu(x, \infty) \le \frac{\|\rho\|_\alpha}{\varepsilon^\alpha \rho(\varepsilon, \infty)} x^{-\alpha}, \qquad x > 0.$$

Hence, the function $h(x) = x^\alpha \nu(x, \infty)$, $x > 0$, is bounded, nonnegative and right-continuous. In terms of this function, (2.7) can be written as

$$\int_0^\infty h(xy^{-1}) y^\alpha \rho(dy) = \|\rho\|_\alpha \qquad \text{for every } x > 0. \tag{2.8}$$

Next, we transform the functional equation (2.8) into additive form by a change of variable. Define a $\sigma$-finite measure on $\mathbb{R}$ by $\mu(dx) = e^{\alpha x}(\rho \circ \log^{-1})(dx)$ and a bounded nonnegative measurable function by $f(x) = h(e^x)$, $x \in \mathbb{R}$. Then (2.8) is transformed into the equivalent equation

$$\int_{-\infty}^\infty f(x - y) \mu(dy) = \|\rho\|_\alpha, \qquad x \in \mathbb{R} \tag{2.9}$$

and upon defining a bounded measurable function by $f_0(x) = f(x) - 1$, we obtain the functional equation

$$(f_0 * \mu)(x) := \int_{-\infty}^\infty f_0(x - y) \mu(dy) = 0 \qquad \text{for all } x \in \mathbb{R}. \tag{2.10}$$

Next, we introduce some basic notions of the theory of generalized functions. We give [32] as a general reference on this theory where one can also find the notions introduced below. Denote by $\mathcal{D}$ the space of $C^\infty$ complex-valued functions with compact support. This space is a subset of the space $\mathcal{S}$ containing the $C^\infty$ complex-valued functions of rapid descent. Let $\mathcal{S}'$ be the dual to $\mathcal{S}$, the space of distributions of slow growth (or tempered distributions). We will use the notation $\langle g, \phi \rangle$ for the action of a tempered distribution $g$ on a test function $\phi \in \mathcal{S}$ and $\hat{g} \in \mathcal{S}'$ will stand for the Fourier transform of $g \in \mathcal{S}'$. We view $\mu$, $f_0$ and $f_0 * \mu$ as elements of $\mathcal{S}'$. Further, if $\mu_n(dx) = \mathbf{1}_{[-n,n]}(x)\mu(dx)$ is the restriction of $\mu$ to the interval $[-n, n]$, $n \geq 1$, then we can view each $\mu_n$ as an element of $\mathcal{S}'$ as well. Moreover, each $\mu_n$ is a distribution with compact support.

For $\phi \in \mathcal{D} \subset \mathcal{S}$, its Fourier transform $\hat{\phi}$ is in $\mathcal{S}$ as well. Now an application of the Lebesgue dominated convergence theorem gives us

$$\lim_{n \to \infty} \langle f_0 * \mu_n, \hat{\phi} \rangle = \langle f_0 * \mu, \hat{\phi} \rangle = 0.$$

Notice that for every $n \geq 1$,

$$\langle f_0 * \mu_n, \hat{\phi} \rangle = \langle \widehat{f_0 * \mu_n}, \phi \rangle = \langle \hat{f_0} \hat{\mu}_n, \phi \rangle = \langle \hat{f_0}, \hat{\mu}_n \phi \rangle,$$

where the first equality is the Parseval identity, the second one follows from Theorem 7.9-1 in [32], and the last identity is the statement of the fact that the Fourier transform of a distribution with compact support is an infinitely smooth function. By assumption (2.3), we get

$$\int_{-\infty}^{\infty} e^{\delta|x|} \mu(dx) = \int_0^1 y^{\alpha-\delta} \rho(dy) + \int_1^\infty y^{\alpha+\delta} \rho(dy) < \infty.$$

Hence, all moments of $\mu$ are finite and its Fourier transform $\hat{\mu}$ itself is in $C^\infty$. Thus, $\hat{\mu}_n \phi \to \hat{\mu} \phi$ in $\mathcal{D}$ which implies the convergence in $\mathcal{S}$. Since $\hat{f_0} \in \mathcal{S}'$, we conclude from (2.10) that for every $\phi \in \mathcal{D}$,

$$(2.11) \qquad \langle \hat{f_0}, \hat{\mu} \phi \rangle = \lim_{n \to \infty} \langle \hat{f_0}, \hat{\mu}_n \phi \rangle = 0.$$

Assumption (2.4) implies that $\phi = \psi / \hat{\mu} \in \mathcal{D}$ provided $\psi \in \mathcal{D}$. Substitution into (2.11) yields

$$(2.12) \qquad \langle \hat{f_0}, \psi \rangle = 0 \qquad \text{for all } \psi \in \mathcal{D}.$$

Since $\mathcal{D}$ is dense in $\mathcal{S}$, relation (2.12) extends to the whole $\mathcal{S}$. Therefore, $\hat{f_0}$ is the zero tempered distribution, and by the uniqueness of the Fourier transform we conclude that so is $f_0$. In terms of functions, $f_0(x) = 0$ a.e. Since $f_0$ is right-continuous, we see that $f_0(x) = 0$ for all $x$. Unwrapping the definition of $f_0$, we obtain $x^\alpha \nu(x, \infty) = 1$ for all $x > 0$, implying that $\nu = \nu_\alpha$, and so the measure $\rho$ has the cancellation property. This proves the first part of the theorem.

Now suppose that the left-hand side of (2.4) vanishes for some $\theta_0 \in \mathbb{R}$. Hence,

$$\int_0^\infty (xy^{-1})^{i\theta_0} y^\alpha \rho(dy) = 0, \qquad x > 0.$$

Taking the real and imaginary parts of this integral, we see that

$$\int_0^\infty g(xy^{-1}) y^\alpha \rho(dy) = \|\rho\|_\alpha, \qquad x > 0,$$

where the function $g$ is defined by (2.6). Consequently, if $\nu$ is given by (2.5), then for every $0 < u < v$,

$$
\begin{aligned}
(\nu \circledast \rho)(u, v] &= \int_0^\infty \left( \int_{uy^{-1}}^{vy^{-1}} g(z) \alpha z^{-(\alpha+1)} \, dz \right) \rho(dy) \\
&= \int_u^v \left( \int_0^\infty g(xy^{-1}) y^\alpha \rho(dy) \right) \alpha x^{-(\alpha+1)} \, dx = \|\rho\|_\alpha \nu_\alpha(u, v],
\end{aligned}
$$

showing that $\nu \circledast \rho = \nu_\alpha \circledast \rho$. Thus, $\rho$ lacks the cancellation property.

Next, we consider the case $\alpha = 0$. We follow the proof for $\alpha > 0$ with some modifications. Assumption (2.3) implies that $\rho$ must be a finite measure with total mass $\|\rho\|_0$. Suppose that (2.4) holds. Define the function

$$
h(x) = x \int_{(x, \infty)} y^{-1} \nu(dy), \qquad x > 0.
$$

Equation (2.7) gives the identity

$$
\int_0^\infty h(xy^{-1}) \rho(dy) = \|\rho\|_0, \qquad x > 0,
$$

from which we deduce that $h$ is bounded. Indeed, if $\rho(\varepsilon, \varepsilon^{-1}) > 0$, then

$$
h(x) \le \frac{\|\rho\|_0}{\varepsilon^2 \rho(\varepsilon, \varepsilon^{-1})}, \qquad x > 0.
$$

The rest of the proof is the same as in the case $\alpha > 0$.

The case $\alpha < 0$ can be immediately reduced to the already solved case $\alpha > 0$ by applying the transformation $x \mapsto x^{-1}$ to measures and their multiplicative convolutions on $(0, \infty)$.  $\square$

The essence of the analysis of the cancellation property of a $\sigma$-finite measure $\rho$ in the proof of Theorem 2.1 above is as to whether the functional equation (2.9) [or its equivalent form (2.10)] have nonconstant solutions. Upon defining a measure of total mass 2 by $\nu = \|\rho\|_\alpha^{-1} \mu + \delta_{\{0\}}$, the question can be stated in the equivalent form: when does the functional equation

$$
(2.13) \qquad \int_{-\infty}^\infty f(x - y) \nu(dy) = f(x) + 1 \qquad \text{for all } x \in \mathbb{R},
$$

have a nonnegative bounded right-continuous solution $f$ different from $f \equiv 1$ a.e.? The functional equation (2.13) is close to the integrated Cauchy functional equation (2.1.1) of [24], apart from the nonhomogeneity introduced by the second term in the right-hand side above. Furthermore, in integrated Cauchy functional equations, one often assumes $\nu(\{0\}) < 1$. Despite the similarity, (2.13) seems to be outside of the framework of existing theory on integrated Cauchy and related functional equations; the stability theorem, Theorem 4.3.1 in [24], for example, requires $\nu$ to have total mass not exceeding 1.

Explicit sufficient conditions for the cancellation property of $\rho$ are available in the presence of atoms.

COROLLARY 2.2. *Let $\alpha \in \mathbb{R}$ and $\rho$ be a nonzero $\sigma$-finite measure satisfying* (2.3) *with an atom, that is, $w_0 = \rho(\{x_0\}) > 0$ for some $x_0 > 0$.*

(i) *If the condition*

$$\int_{y \neq x_0} y^\alpha \rho(dy) < w_0 x_0^\alpha \tag{2.14}$$

*holds, then the measure $\rho$ has the cancellation property with respect to $\mathcal{N} = \{\nu_\alpha\}$.*

(ii) *Suppose that*

$$\int_{y \neq x_0} y^\alpha \rho(dy) = w_0 x_0^\alpha. \tag{2.15}$$

*Then the measure $\rho$ lacks the cancellation property with respect to $\mathcal{N} = \{\nu_\alpha\}$ if and only if for some $\theta_0 > 0$, $\rho(S_{x_0}^c) = 0$, where $S_{x_0}^c$ is the complement of the set*

$$S_{x_0} = \{x_0\} \cup \{x_0 e^{\pi(2k+1)/\theta_0}, k \in \mathbb{Z}\}. \tag{2.16}$$

PROOF. (i) From (2.14) and the elementary inequality

$$\left| \int_0^\infty y^{\alpha+i\theta} \rho(dy) \right| \geq |w_0 x_0^{\alpha+i\theta}| - \left| \int_{y \neq x_0} y^{\alpha+i\theta} \rho(dy) \right|$$

$$\geq w_0 x_0^\alpha - \int_{y \neq x_0} y^\alpha \rho(dy) > 0,$$

we conclude that (2.4) holds. The cancellation property of $\rho$ follows from Theorem 2.1.

(ii) Suppose that the left-hand side of (2.4) vanishes for some $\theta_0$, which gives

$$\int_{y \neq x_0} y^{\alpha+i\theta_0} \rho(dy) = -w_0 x_0^{\alpha+i\theta_0}. \tag{2.17}$$

We may always take $\theta_0 > 0$. By (2.15), we get

$$\left| \int_{y \neq x_0} y^{\alpha+i\theta_0} \rho(dy) \right| = \int_{y \neq x_0} |y^{\alpha+i\theta_0}| \rho(dy),$$

which by Theorem 2.1.4 in [19] says that $y^{i\theta_0}$ must be a.s. a constant with respect to the measure $\rho$ on $(0, \infty) \setminus \{x_0\}$. By (2.17), this constant must be equal to $-x_0^{i\theta_0}$, and so the measure $\rho$ cannot have mass outside of the set $S_{x_0}$ in (2.16). Conversely, if the measure $\rho$ is concentrated on the set $S_{x_0}$, then (2.17) clearly holds, and so the left-hand side of (2.4) vanishes at the point $\theta_0$. Now appeal to Theorem 2.1. □

The following result relies on Theorem 2.1. It is fundamental in our studies of the inverse problems for regular variation given in the subsequent sections. Interestingly, there does not seem to be an analogous statement in the case $\alpha = 0$.

THEOREM 2.3.  *Let $\rho$ be a nonzero $\sigma$-finite measure satisfying* (2.3) *for some $\alpha > 0$ and assume that condition* (2.4) *is satisfied. Suppose that, for some $\sigma$-finite measure $\nu$ on $(0, \infty)$, the measure $\nu \circledast \rho$ in* (2.1) *has a regularly varying tail with exponent $\alpha$ and*

$$(2.18) \qquad \lim_{b \to 0} \limsup_{x \to \infty} \frac{\int_0^b \rho(x/y, \infty) \nu(dy)}{(\nu \circledast \rho)(x, \infty)} = 0.$$

*Then the measure $\nu$ has regularly varying tail with exponent $\alpha$ as well, and*

$$(2.19) \qquad \lim_{x \to \infty} \frac{(\nu \circledast \rho)(x, \infty)}{\nu(x, \infty)} = \int_0^\infty y^\alpha \rho(dy).$$

*Conversely, if* (2.4) *does not hold, then there exists a $\sigma$-finite measure $\nu$ on $(0, \infty)$ without a regularly varying tail, such that the measure $\nu \circledast \rho$ has regularly varying tail with exponent $\alpha$ and* (2.18) *holds.*

PROOF.  Suppose that (2.4) is satisfied and let $\nu$ be a measure as described in the theorem. For $n \geq 1$, define

$$a_n = \min\{(\nu \circledast \rho)(n, \infty), 1\},$$

and assume without loss of generality that $\max\{n : a_n = 1\} = 1$. Moreover, define the sequence of Radon measures $(\nu_n)$ on $(0, \infty]$ by

$$\nu_n(x, \infty] = \nu_n(x, \infty) = a_n^{-1} \nu(nx, \infty), \qquad x > 0.$$

Then for every $x > 0$,

$$(2.20) \qquad \int_0^\infty \nu_n(x/y, \infty) \rho(dy) = \frac{(\nu \circledast \rho)(nx, \infty)}{(\nu \circledast \rho)(n, \infty)} \to x^{-\alpha} \qquad \text{as } n \to \infty.$$

Fix $y_0 > 0$ such that $C_0 := \rho[y_0, \infty) > 0$. Then for any $x > 0$ and $n \geq 1$,

$$(2.21) \qquad \nu_n(x, \infty) \leq C_0^{-1} \frac{(\nu \circledast \rho)(nxy_0, \infty)}{(\nu \circledast \rho)(n, \infty)}.$$

In particular,

$$\limsup_{n \to \infty} \nu_n(x, \infty) \leq [C_0 y_0^\alpha]^{-1} x^{-\alpha}.$$

By Proposition 3.16 in [25], this implies that the sequence $(\nu_n)$ is vaguely relatively compact. We claim that this sequence converges vaguely. To this end, it is enough to show that all vague subsequential limits $\gamma_*$ of the sequence $(\nu_n)$ coincide.

Assume that $\gamma_*$ is the vague limit of $(\nu_{n_k})$ for some integer subsequence $(n_k)$. For $M > 0$, we write

$$\int_0^\infty \nu_n(x/y, \infty) \rho(dy) = \left( \int_{y \leq Mn} + \int_{y > Mn} \right) \nu_n(x/y, \infty) \rho(dy).$$

Let $0 < \delta < \alpha$ be such that (2.3) holds. By (2.21), we have that for $x > 0$, $n$ and $y \leq Mn$,

$$\nu_n(x/y, \infty) \leq C_0^{-1} \frac{(\nu \circledast \rho)(nxy_0/y, \infty)}{(\nu \circledast \rho)(n, \infty)} \leq C(x, y_0)(y^{\alpha - \delta} \vee y^{\alpha + \delta}),$$

where the last inequality follows from the Potter bounds, compare [4], Theorem 1.5.6 or Proposition 0.8(ii) in [25]. Using (2.3) and the Lebesgue dominated convergence theorem along the subsequence $(n_k)$, we obtain

$$(2.22) \qquad \int_{y \leq Mn_k} \nu_{n_k}(x/y, \infty) \rho(dy) \to \int_0^\infty \gamma_*(x/y, \infty) \rho(dy)$$

provided $\rho(D_x) = 0$, where $D_x = \{y > 0 : \gamma_*(\{x/y\}) > 0\}$. Since $D_x$ is at most countable, $\rho(D_x) > 0$ if and only if $D_x$ contains an atom of $\rho$. Since the set of $x \in \mathbb{R}$ such that $D_x$ contains a specific atom of $\rho$ is at most countable, we have $\rho(D_x) = 0$ with the possible exception of a countable set of $x > 0$. Therefore, (2.22) is valid for all, but countably many $x > 0$.

On the other hand, by assumption (2.3), we get for $M \geq 1$,

$$\int_{y > Mn} \nu_n(x/y, \infty) \rho(dy)$$
$$= \frac{\rho(Mn, \infty)}{a_n} \nu(M^{-1}x, \infty) + \frac{1}{a_n} \int_0^{x/M} \rho(nx/y, \infty) \nu(dy)$$
$$\leq \frac{\int_1^\infty y^{\alpha + \delta} \rho(dy)}{(Mn)^{\alpha + \delta} a_n} \nu(M^{-1}x, \infty) + \frac{1}{a_n} \int_0^{x/M} \rho(nx/y, \infty) \nu(dy).$$

Using Proposition 0.8(ii) in [25], we get

$$\limsup_{n \to \infty} \int_{y > Mn} \nu_n(x/y, \infty) \rho(dy) \leq x^{-\alpha} \limsup_{u \to \infty} \frac{\int_0^{x/M} \rho(u/y, \infty) \nu(dy)}{(\nu \circledast \rho)(u, \infty)}.$$

Combining this with (2.18) and (2.22), we conclude that

$$\lim_{k \to \infty} \int_0^\infty \nu_{n_k}(x/y, \infty) \rho(dy) = \int_0^\infty \gamma_*(x/y, \infty) \rho(dy)$$

for all but countably many $x > 0$. The limit (2.20) now says that

$$(2.23) \qquad \int_0^\infty \gamma_*(x/y, \infty) \rho(dy) = x^{-\alpha} = \nu_\alpha(x, \infty)$$

for all but countably many $x > 0$, but by right continuity, (2.23) extends to all $x > 0$. This equation can also be written as

$$\|\rho\|_\alpha (\gamma_* \circledast \rho) = \nu_\alpha \circledast \rho,$$

see (2.7), where as usual, $\|\rho\|_\alpha = \int_0^\infty y^\alpha \rho(dy)$. Applying Theorem 2.1, we get that $\gamma_* = \|\rho\|_\alpha^{-1} \nu_\alpha$. Since this is true for any vague subsequential limit of the sequence $(\nu_n)$, we conclude that there is, in fact, full vague convergence $\nu_n \xrightarrow{v} \|\rho\|_\alpha^{-\alpha} \nu_\alpha$. This implies both regular variation of the tail of the measure $\nu$ and the relation (2.19), thus completing the proof of the first part of the theorem.

In the opposite direction, suppose that the left-hand side of (2.4) vanishes for some $\theta_0 \in \mathbb{R}$. We may assume that $\theta_0 > 0$. Take real $a, b$ such that $0 < a^2 + b^2 < 1$. Define a $\sigma$-finite measure $\nu$ by (2.5) with the function $g$ given by (2.6). By Theorem 2.1, $\nu \circledast \rho = \|\rho\|_\alpha \nu_\alpha$, which clearly, has a regularly varying tail with exponent $\alpha$. Furthermore, (2.18) holds as well. We claim that the measure $\nu$ does not have a regularly varying tail.

Indeed, suppose on the contrary that the tail of the measure $\nu$ is regularly varying with some index $-\beta$. Note that for every $x > 0$, $g(rx) = g(x)$, where $r = e^{2\pi/\theta_0} > 1$. This implies $\nu(r^m x, \infty) = r^{-\alpha m} \nu(x, \infty)$ for all $m = 1, 2, \ldots$ and $x > 0$. It follows that $\beta = \alpha$. Since

$$\lim_{m \to \infty} \frac{\nu(r^m x, \infty)}{\nu(r^m, \infty)} = x^{-\alpha},$$

$\nu(x, \infty) = x^{-\alpha} \nu(1, \infty)$ for all $x > 0$. This contradicts the definition of $\nu$ in (2.5). Therefore, the measure $\nu$ does not have a regularly varying tail. $\square$

REMARK 2.4. Condition (2.18) says that the contribution of the left tail of $\nu$ to the right tail of the product convolution is negligible. This condition is automatic in many situations of interest. It is trivially satisfied if the measure $\rho$ is supported on a bounded interval $(0, B]$ for some finite $B$. The condition always holds when the measure $\nu$ is finite in a neighborhood of the origin. More generally, if $\int_0^1 y^{\alpha+\delta} \nu(dy) < \infty$, with $\delta \in (0, \alpha)$ satisfying (2.3), then (2.18) holds. Then we have for $b \in (0, 1)$ and $x > 1$,

$$\int_0^b \rho(x/y, \infty) \nu(dy) \leq x^{-(\alpha+\delta)} \int_1^\infty y^{\alpha+\delta} \rho(dy) \int_0^1 y^{\alpha+\delta} \nu(dy) = C x^{-(\alpha+\delta)},$$

which implies (2.18) because of the assumed regular variation of the measure $\nu \circledast \rho$.

**3. The inverse problem for weighted sums.** In this section, we consider the weighted sums of Example 1.1, and address the question, to what extent regular variation of the sum in (1.2) implies regular variation of the noise variables. The following is the main result of this section.

THEOREM 3.1. *Assume that $\alpha > 0$ and $(Z_j)$ is a sequence of i.i.d. random variables.*

(i) *Suppose that $(\psi_j)$ is a sequence of positive coefficients satisfying*

$$(3.1) \qquad \sum_{j=1}^{\infty} \psi_j^{\alpha-\delta} < \infty \qquad \text{for some } 0 < \delta < \alpha.$$

*If $\sum_{j=1}^{\infty} \psi_j = \infty$, assume additionally that*

$$(3.2) \qquad \limsup_{x \to \infty} \frac{P(Z < -x)}{P(Z > x)} < \infty.$$

*Assume that the series $X = \sum_{j=1}^{\infty} \psi_j Z_j$ converges a.s., and that $X$ is regularly varying with exponent $\alpha$. If*

$$(3.3) \qquad \sum_{j=1}^{\infty} \psi_j^{\alpha+i\theta} \neq 0 \qquad \text{for all } \theta \in \mathbb{R},$$

*then a generic noise variable $Z$ is regularly varying with exponent $\alpha$ as well, and (1.3) holds.*

(ii) *Suppose that $(\psi_j)$ is a sequence of positive coefficients satisfying (3.1). If (3.3) fails to hold, then there exists a random variable $Z$ that is not regularly varying, the series $X = \sum_{j=1}^{\infty} \psi_j Z_j$ converges a.s., and $X$ is regularly varying with exponent $\alpha$.*

PROOF. (i) We start with the observation that for some $C > 0$

$$(3.4) \qquad P(Z > x) \leq C P(X > x) \qquad \text{for all } x > 0.$$

This bound follows from the inequality

$$P(X > x) \geq P(\psi_1 Z_1 > x + M) P\left( \sum_{j=2}^{\infty} \psi_j Z_j \geq -M \right)$$

and the regular variation of $X$, where $M > 0$ is such that the right-hand side is positive.

Our next goal is to show that condition (3.2) may be assumed to hold under the conditions of the theorem. When $\sum_{j=1}^{\infty} \psi_j = \infty$, we have already assumed this condition.

If $\alpha \leq 1$, then (3.4), assumption (3.1) and the 3-series theorem (see Lemma 6.9 and Theorem 6.1 in [22]) guarantee that

$$(3.5) \qquad X^{(+)} := \sum_{j=1}^{\infty} \psi_j(Z_j)_+ < \infty \qquad \text{a.s.;}$$

see, for example, Section 4.5 in [25]. Here, $a_{\pm} = \max(\pm a, 0)$. Therefore, the series $X^{(-)} := \sum_{j=1}^{\infty} \psi_j(Z_j)_-$ also converges a.s. Clearly, for any $M > 0$ and $x > M$, we have

$$(3.6) \qquad P(X^{(+)} > x) \geq P(X > x) \geq P(X^{(+)} > x + M, X^{(-)} \leq M).$$

Notice that the random variables $X^{(+)}$ and $-X^{(-)}$ are associated, as both are nondecreasing functions of i.i.d. random variables $Z_j$, $j = 1, 2, \ldots$; see [16]. Therefore,

$$
\begin{aligned}
(3.7) \qquad P(X^{(+)} > x + M, X^{(-)} \leq M) &= P(X^{(+)} > x + M, -X^{(-)} \geq -M) \\
&\geq P(X^{(+)} > x + M)P(X^{(-)} \leq M).
\end{aligned}
$$

It follows from (3.6), (3.7) and the regular variation of the tail of $X$ that

$$1 \geq \limsup_{x \to \infty} \frac{P(X > x)}{P(X^{(+)} > x)} \geq \liminf_{x \to \infty} \frac{P(X > x)}{P(X^{(+)} > x)} \geq P(X^{(-)} \leq M).$$

Letting $M \to \infty$, we conclude that $P(X > x) \sim P(X^{(+)} > x)$ as $x \to \infty$. In particular, we may assume, without loss of generality, that (3.2) holds in this case.

If $\alpha > 1$ and $\sum_{j=1}^{\infty} \psi_j < \infty$, then (3.4) still implies (3.5) (just take the expectation), and the above argument allows us to assume that (3.2) holds in this case as well.

For the reasons given, throughout the proof of the first part of the theorem, we may and will assume (3.2). Together with (3.4), this shows that for some $C > 0$

$$P(|Z| > x) \leq CP(X > x) \qquad \text{for all } x > 0.$$

For $K \geq 1$, let $X^{(K)} = \sum_{j=K+1}^{\infty} \psi_j Z_j$. We next claim that

$$(3.8) \qquad \lim_{K \to \infty} \limsup_{x \to \infty} \frac{P(|X^{(K)}| > x)}{P(X > x)} = 0.$$

The proof is identical to the argument for (A.5) in [21] and, therefore, omitted.

We claim that under the assumptions of part (i),

$$(3.9) \qquad P(X > x) \sim \sum_{j=1}^{\infty} P(Z > x/\psi_j) \qquad \text{as } x \to \infty.$$

Because of (3.8), it is enough to prove (3.9) under the assumption that the coefficients $\psi_j$ vanish after, say, $j = K$. In that case, by a Bonferroni argument, for $\varepsilon > 0$,

$$P(X > x) \geq P\left(\bigcup_{j=1}^{K} \{Z_j > (1+\varepsilon)x/\psi_j, Z_i \leq \varepsilon x/(K\psi_i) \text{ for all } i \neq j\}\right)$$

$$\geq \sum_{j=1}^{K} P(Z_j > (1+\varepsilon)x/\psi_j, Z_i \leq \varepsilon x/(K\psi_i) \text{ for all } i \neq j)$$

$$\qquad - \left(\sum_{j=1}^{K} P(Z_j > (1+\varepsilon)x/\psi_j)\right)^2$$

$$\geq \sum_{j=1}^{K} \Bigg[ P(Z_j > (1+\varepsilon)x/\psi_j)$$

$$\qquad\qquad - \sum_{i=1,\ldots,K, i \neq j} P(Z_j > \varepsilon x/(K\psi_j), Z_i > \varepsilon x/(K\psi_i)) \Bigg]$$

$$\qquad - \left(\sum_{j=1}^{K} P(Z_j > \varepsilon x/(K\psi_j))\right)^2$$

$$\geq \sum_{j=1}^{K} P(Z_j > (1+\varepsilon)x/\psi_j) - 2\left(\sum_{j=1}^{K} P(Z_j > \varepsilon x/(K\psi_j))\right)^2.$$

Using the upper bound (3.4) and the arbitrariness of $\varepsilon > 0$ gives us the bound

$$\liminf_{x \to \infty} \frac{P(X > x)}{\sum_{j=1}^{K} P(Z > x/\psi_j)} \geq 1.$$

Similarly, in the other direction: for $0 < \varepsilon < 1$,

$$P(X > x) \leq P\Bigg[\bigcup_{j=1}^{K} \{Z_j > (1-\varepsilon)x/\psi_j\}$$

$$\qquad \cup \bigcup_{i,j=1,\ldots,K, i \neq j} \{Z_i > \varepsilon x/(K\psi_i), Z_j > \varepsilon x/(K\psi_j)\}\Bigg]$$

$$\leq \sum_{j=1}^{K} P(Z_j > (1-\varepsilon)x/\psi_j) + \left(\sum_{j=1}^{K} P(Z_j > \varepsilon x/(K\psi_j))\right)^2.$$

Together with the upper bound (3.4) and the arbitrariness of $\varepsilon > 0$, we obtain the second bound

$$\limsup_{x \to \infty} \frac{P(X > x)}{\sum_{j=1}^{K} P(Z > x/\psi_j)} \leq 1,$$

which proves (3.9).

Define the measures

(3.10)
$$\nu(\cdot) = P(Z \in \cdot) \quad \text{and}$$

$$\rho(B) = \sum_{j=1}^{\infty} \mathbf{1}_B(\psi_j) \qquad \text{for a Borel subset } B \text{ of } (0, \infty).$$

The just proved relation (3.9) implies that the measure $\nu \circledast \rho$ has a regularly varying tail with exponent $\alpha$. Assumption (3.1) implies (2.3), and condition (2.18) also holds; see Remark 2.4. Therefore, Theorem 2.3 applies and gives us both the regular variation of the tail of $Z$, and the tail equivalence (1.3). This completes the proof of part (i).

(ii) Assuming that (3.3) fails, by the construction in Theorem 2.3, we can find a $\sigma$-finite measure $\nu$ that does not have a regularly varying tail, such that

$$\sum_{j=1}^{\infty} \nu(x/\psi_j, \infty) = cx^{-\alpha} \qquad \text{for all } x > 0, \text{ some } c > 0.$$

Choose $b > 0$ large enough so that $\nu(b, \infty) \leq 1$, and define a probability law on $(0, \infty)$ by

$$\mu(B) = \nu(B \cap (b, \infty)) + [1 - \nu(b, \infty)]\mathbf{1}_B(1) \qquad \text{for a Borel set } B.$$

We then have

(3.11)
$$\sum_{j=1}^{\infty} \mu(x/\psi_j, \infty) = cx^{-\alpha} \qquad \text{for all } x \text{ large enough.}$$

Let $Z$ have the law

$$P(Z \in A) = \tfrac{1}{2}\mu(A \cap (0, \infty)) + \tfrac{1}{2}\mu(-A \cap (0, \infty))$$

for Borel subsets $A$ of the reals, and let $Z_j$, $j = 1, 2, \ldots$, be i.i.d. copies of $Z$. Then $(Z_j)$ is a sequence of i.i.d. symmetric random variables, whose tails satisfy

$$P(|Z| > x) \leq cx^{-\alpha} \qquad \text{for } x > 0,$$

where $c$ is a positive constant. By assumption (3.1) and Lemma A.3 in [21], the series $X = \sum_{j=1}^{\infty} \psi_j Z_j$ converges.

We can employ the argument used in part (i) of the theorem to verify that as $x \to \infty$,

$$P(X > x) \sim \sum_{j=1}^{\infty} P(Z > x/\psi_j) = \frac{1}{2} \sum_{j=1}^{\infty} \mu(x/\psi_j, \infty) = \frac{c}{2} x^{-\alpha}$$

by (3.11). Therefore, $X$ has a regularly varying tail with exponent $\alpha$, thus giving us the required construction. $\square$

EXAMPLE 3.2. A stationary AR(1) process with a positive coefficient has the one-dimensional marginal distribution given in the form of the series (1.2) with $\psi_j = \beta^{j-1}, j = 1, 2, \ldots$, for some $0 < \beta < 1$. For these coefficients, the sum in the right-hand side of (3.3) becomes

$$\sum_{j=0}^{\infty} \beta^{j(\alpha + i\theta)} = \frac{1}{1 - \beta^{\alpha + i\theta}},$$

which does not vanish. Therefore, regular variation of the marginal distribution of such a process is equivalent to the regular variation of the noise variables.

The perhaps most illuminating application of Theorem 3.1 concerns finite sums. In this case, all the issues related to the convergence of infinite sums are eliminated. We consider such finite sums in the remainder of this section.

We say that a set of $q \geq 2$ positive coefficients $\psi_1, \ldots, \psi_q$ is $\alpha$-*regular variation determining* if any set of i.i.d. random variables $Z_1, \ldots, Z_q$ are regularly varying with exponent $\alpha$ if and only if

$$(3.12) \qquad\qquad X_q = \sum_{j=1}^{q} \psi_j Z_j$$

is regularly varying with exponent $\alpha$. The corresponding notion in the slowly varying case, that is, when $\alpha = 0$, is not of interest: any set of positive coefficients $\psi_1, \ldots, \psi_q$ is 0-regular variation determining. Indeed, suppose that $\psi_1 \geq \psi_2 \geq \cdots \geq \psi_q$ and without loss of generality $\psi_1 = 1$. Then for any $x > 0$,

$$P(X_q > qx) \leq \sum_{j=1}^{q} P(\psi_j Z_j > x) \leq q P(Z > x)$$

and

$$P(X_q > (\psi_q/2)x) \geq P\left( \bigcup_{j=1}^{q} \left\{ \psi_j Z_j > \psi_q x, \sum_{i \neq j} |\psi_i Z_i| \leq (\psi_q/2)x \right\} \right)$$

$$= \sum_{j=1}^{q} P(\psi_j Z_1 > \psi_q x) P\left(\sum_{i \neq j} |\psi_i Z_i| \leq (\psi_q/2)x\right)$$

$$\geq q P(Z > x) P\left(\sum_{j=1}^{q-1} |Z_j| \leq (\psi_q/2)x\right).$$

Hence,

$$\frac{1}{q} P(X_q > qx) \leq P(Z > x) \leq \frac{1}{q} \frac{P(X_q > (\psi_q/2)x)}{P(\sum_{j=1}^{q-1} |Z_j| \leq (\psi_q/2)x)}$$

and the slow variation of the tail of $X_q$ implies that $P(Z > x) \sim q^{-1}P(X_q > x)$ as $x \to \infty$. Therefore, we always assume $\alpha > 0$.

Next, we reformulate Theorem 3.1 and Corollary 2.2 for finite sums.

THEOREM 3.3. (i) *A set of positive coefficients $\psi_1, \ldots, \psi_q$ is $\alpha$-regular variation determining if and only if*

$$(3.13) \qquad \sum_{j=1}^{q} \psi_j^{\alpha+i\theta} \neq 0 \qquad \text{for all } \theta \in \mathbb{R}.$$

(ii) *Suppose that $\psi_{j_1} = \cdots = \psi_{j_k}$ for some subset $\{j_1, \ldots, j_k\}$ of $\{1, \ldots, q\}$, and $\psi_j \neq \psi_{j_1}$ if $j \in \{1, \ldots, q\} \setminus \{j_1, \ldots, j_k\}$. If*

$$(3.14) \qquad \sum_{\psi_j \neq \psi_{j_1}} \psi_j^{\alpha} < k\psi_{j_1}^{\alpha},$$

*then $\psi_1, \ldots, \psi_q$ are $\alpha$-regular variation determining. If, on the other hand,*

$$(3.15) \qquad \sum_{\psi_j \neq \psi_{j_1}} \psi_j^{\alpha} = k\psi_{j_1}^{\alpha},$$

*then $\psi_1, \ldots, \psi_q$ are $\alpha$-regular variation determining if and only if there is no $\theta_0 > 0$ such that for every $j \in \{1, \ldots, q\}$ with $\psi_j \neq \psi_{j_1}$,*

$$(3.16) \qquad \frac{\log \psi_j - \log \psi_{j_1}}{\pi \theta_0} \quad \text{is an odd integer.}$$

REMARK 3.4. An immediate conclusion from part (i) of Theorem 3.3 is that the coefficients $\psi_1, \ldots, \psi_q$ are $\alpha$-regular variation determining if and only if the coefficients $\psi_1^{\alpha}, \ldots, \psi_q^{\alpha}$ are 1-regular variation determining.

EXAMPLE 3.5. The case $\psi_1 = \cdots = \psi_q$ is a well known $\alpha$-regular variation determining one. It corresponds to the so-called *convolution root closure property* of subexponential distributions; see [14]. Distributions on $(0, \infty)$ with a regularly varying right tail constitute a subclass of the subexponential distributions. In a sense, the property of $\alpha$-regular variation determination can be understood as a natural extension of the convolution root closure property to weighted sums.

EXAMPLE 3.6.   For $q = 2$, one can always apply part (ii) of Theorem 3.3. If $\psi_1 \neq \psi_2$, then (3.14) holds with $k = 1$ and $\psi_{j_1} = \max(\psi_1, \psi_2)$. If $\psi_1 = \psi_2$, then (3.14) holds with $k = 2$. Therefore, any set of 2 positive coefficients is $\alpha$-regular variation determining.

EXAMPLE 3.7.   Take an arbitrary $q \geq 2$, and suppose that for some $1 \leq m \leq q$, $\psi_1 = \cdots = \psi_m \neq \psi_{m+1} = \cdots = \psi_q$. Once again, part (ii) of Theorem 3.3 applies. If

$$(3.17) \qquad m\psi_1^\alpha \neq (q - m)\psi_q^\alpha,$$

then (3.14) holds with $k = m$ or $k = q - m$, depending on which side of (3.17) is greater, and the coefficients $\psi_1, \ldots, \psi_q$ are $\alpha$-regular variation determining. On the other hand, if (3.17) fails, then (3.15) holds with $k = m$, and the coefficients $\psi_1, \ldots, \psi_q$ are not $\alpha$-regular variation determining because (3.16) holds with

$$\theta_0 = \frac{|\log \psi_1 - \log \psi_q|}{\pi}.$$

EXAMPLE 3.8.   Let $q = 3$, and assume that $\alpha = 1$ (by Remark 3.4, we can switch at will between different values of exponent of regular variation $\alpha$).

Since the property of $\alpha$-regular variation determination is invariant under multiplication of all coefficients with the same positive number, we assume $\psi_3 = 1$ and $\psi_1, \psi_2 < 1$, $\psi_1 \neq \psi_2$. Otherwise we are, once again, in the situation of Example 3.7. By part (ii) of Theorem 3.3, the coefficients are 1-regular variation determining if

$$\psi_1 + \psi_2 < 1 \quad \text{or if}$$

$$\psi_1 + \psi_2 = 1 \quad \text{and} \quad \frac{\log \psi_1}{\log \psi_2} \neq \frac{2m_1 + 1}{2m_2 + 1} \qquad \text{for some } m_1, m_2 = 0, 1, \ldots.$$

Suppose that $\psi_1 + \psi_2 > 1$. By part (i) of Theorem 3.3, the coefficients fail to be $\alpha$-regular variation determining if and only if for some real $\theta$

$$\psi_1 e^{i\theta \log \psi_1} + \psi_2 e^{i\theta \log \psi_2} = -1,$$

which is equivalent to the system of two equations,

$$(3.18) \qquad \psi_1 \cos(\theta \log \psi_1) + \psi_2 \cos(\theta \log \psi_2) = -1$$

and

$$(3.19) \qquad \psi_1 \sin(\theta \log \psi_1) + \psi_2 \sin(\theta \log \psi_2) = 0.$$

Squaring the two equations and summing them up, one uses the elementary trigonometric formulae to obtain

$$\psi_1^2 + \psi_2^2 + 2\psi_1\psi_2 \cos(\theta \log(\psi_1/\psi_2)) = 1,$$

implying that

$$(3.20) \qquad \theta \log(\psi_1/\psi_2) = \pm \arccos \frac{1 - (\psi_1^2 + \psi_2^2)}{2\psi_1\psi_2} + 2\pi n$$

for some $n \in \mathbb{Z}$. On the other hand, moving the term $\psi_2 \cos(\theta \log \psi_2)$ in (3.18) before squaring it and using (3.19), one obtains the equation

$$\psi_1^2 = 1 + 2\psi_2 \cos(\theta \log \psi_2) + \psi_2^2,$$

so that

$$(3.21) \qquad \theta \log(\psi_2) = \pm \arccos \frac{\psi_1^2 - \psi_2^2 - 1}{2\psi_2} + 2\pi m$$

for some $m \in \mathbb{Z}$. A comparison of (3.20) with (3.21) shows that the coefficients $\psi_1$ and $\psi_2$ must satisfy the relation

$$(3.22) \qquad \frac{\pm \arccos((1 - (\psi_1^2 + \psi_2^2))/(2\psi_1\psi_2)) + 2\pi n}{\log(\psi_1/\psi_2)}$$
$$= \frac{\pm \arccos((\psi_1^2 - \psi_2^2 - 1)/(2\psi_2)) + 2\pi m}{\log(\psi_2)}$$

for some choices of $\pm$ signs and $n, m \in \mathbb{Z}$.

While it appears to be difficult to find explicit solutions to (3.22), it emphasizes that coefficients that do not possess the property of $\alpha$-regular variation determination are "reasonably rare": in the case $q = 3$ they have to lie on a countable set of curves described in (3.22). Some of these curves are presented in Figure 1 in the $(\psi_1, \psi_2)$-coordinates.

## 4. The inverse problem for products of independent random variables.
In this section, we consider products of independent random variables as in Example 1.2. Let $\alpha > 0$ and $Y$ a positive random variable satisfying $EY^{\alpha+\delta} < \infty$ for some $\delta > 0$. We will call $Y$ and its distribution $\alpha$-*regular*

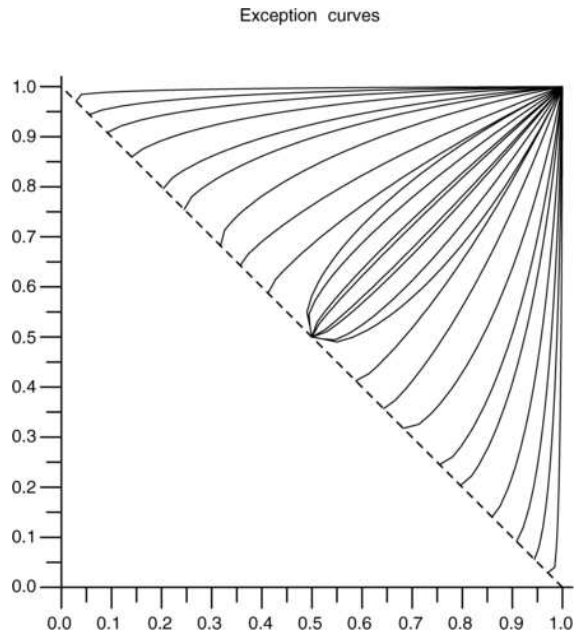

FIG. 1. *Some of the curves satisfying* (3.22).

*variation determining* if the $\alpha$-regular variation of $X = YZ$ for any random variable $Z$ which is independent of $Y$, implies that $Z$ itself has a regularly varying tail with exponent $\alpha$. The corresponding notion in the slowly varying case ($\alpha = 0$) turns out to be not of great interest; we briefly address this point at the end of the section. Until then, we assume that $\alpha > 0$.

The following lemma lists some elementary properties of regular variation determining random variables that follow directly from the definition.

LEMMA 4.1.    (i) *If $Y$ is $\alpha$-regular variation determining, and $\beta > 0$, then $Y^\beta$ is $\alpha/\beta$-regular variation determining.*

(ii) *If $Y_1$ and $Y_2$ are independent $\alpha$-regular variation determining random variables, then so is $Y_1 Y_2$.*

Necessary and sufficient conditions for $\alpha$-regular variation determination of a random variable can be obtained from Theorem 2.3.

THEOREM 4.2.    *A positive random variable $Y$ with $EY^{\alpha+\delta} < \infty$ for some $\delta > 0$ is $\alpha$-regular variation determining if and only if*

$$(4.1) \qquad E[Y^{\alpha+i\theta}] \neq 0 \qquad \text{for all } \theta \in \mathbb{R}.$$

Write $F_Y$ for the distribution of $Y$. If $Y$ and $Z$ are independent positive random variables, the quantity

$$P(X \leq x) = \int_0^\infty P(Z \leq x/y) F_Y(dy), \qquad x > 0$$

is sometimes referred to as *Mellin–Stieltjes convolution*; see, for example, [4]. Theorem 4.2 has some close relation with Tauberian and Mercerian theorems; see [5] and the references therein.

In the literature, several examples were found showing that regular variation of the product $X = YZ$ does not necessarily imply regular variation of the factors $Y$ and $Z$. A counterexample close in spirit to the construction in (2.5) and (2.6) can be found in [20] who dedicated this example to Daren Cline. See also [11]. A counterexample of a different type is in the paper [29]. The latter paper also contains a proof of the "if" part of Theorem 4.2 for a particular family of discrete random variables $Y$.

PROOF OF THEOREM 4.2. The only part that requires a proof is that, if condition (4.1) fails, then the measure providing a counterexample in Theorem 2.3 can be chosen to be a finite measure. We use a construction similar to that in the proof of Theorem 3.1. By Theorem 2.3, there exists a $\sigma$-finite measure $\nu$ that does not have a regularly varying tail, such that for all $x > 0$

$$(4.2) \qquad E\nu(x/Y, \infty) = x^{-\alpha}.$$

Define a finite measure $\nu_1$ by

$$\nu_1(B) = \nu(B \cap (b, \infty)), \qquad B \text{ a Borel set,}$$

where $b > 0$ is large enough so that $\nu(b, \infty) < \infty$. By construction, $\nu_1$ does not have a regularly varying tail. We claim that

$$(4.3) \qquad \lim_{x \to \infty} x^\alpha E\nu_1(x/Y, \infty) = 1.$$

Indeed, by (4.2)

$$E\nu_1(x/Y, \infty) = E[\nu(x/Y, \infty)\mathbf{1}(Y < x/b)] + \nu(b, \infty)P(Y \geq x/b)$$
$$= x^{-\alpha} - E[\nu(x/Y, \infty)\mathbf{1}(Y \geq x/b)] + \nu(b, \infty)P(Y \geq x/b).$$

The bound $\nu(x, \infty) \leq cx^{-\alpha}$ valid for all $x > 0$ for some positive constant $c$ says that

$$E[\nu(x/Y, \infty)\mathbf{1}(Y \geq x/b)] \leq cx^{-\alpha}E[Y^\alpha\mathbf{1}(Y \geq x/b)] = o(x^{-\alpha})$$

since $EY^\alpha < \infty$. For the same reason,

$$P(Y \geq x/b) = o(x^{-\alpha})$$

as well, and (4.3) follows. □

Assume $Y$ is a positive random variable satisfying $EY^\alpha < \infty$. Then we can define the $\alpha$-*conjugate to* $Y$ as a random variable $Y^{*,\alpha}$ with law

$$(4.4) \qquad F_{Y^{*,\alpha}}(dx) = ce^{\alpha x}(F_Y \circ \log^{-1})(dx), \qquad x \in \mathbb{R},$$

where $c$ is a normalizing constant required to make $F_{Y^{*,\alpha}}$ a probability measure on $\mathbb{R}$.

COROLLARY 4.3. (i) *A positive random variable $Y$ with $EY^{\alpha+\delta} < \infty$ for some $\delta > 0$ is $\alpha$-regular variation determining if and only if its $\alpha$-conjugate random variable has a nonvanishing characteristic function.*

(ii) *If $\log Y$ is infinitely divisible, then $Y$ is $\alpha$-regular variation determining.*

PROOF. The first part of the corollary is a restatement of Theorem 4.2. For the second part, observe that infinite divisibility of $\log Y$ implies infinite divisibility of the $\alpha$-conjugate random variable of $Y$. Thus, the statement follows from part (i) because the characteristic function of an infinitely divisible random variable does not vanish; see, for example, [19] or [28]. □

As the following examples show, many of the standard distributions are $\alpha$-regular variation determining.

EXAMPLE 4.4 (Gamma and normal random variables). For a $\Gamma(\beta, \lambda)$ random variable,

$$E[Y^{\alpha+i\theta}] = c\frac{\Gamma(\alpha + \beta + i\theta)}{\Gamma(\beta)} \neq 0$$

since the gamma function does not vanish whenever it is defined (here $c$ is a nonzero complex constant). Therefore, a gamma random variable is $\alpha$-regular variation determining for any $\alpha > 0$. Since the square of a centered normal random variable has a gamma distribution, this fact and an appeal to part (i) of Lemma 4.1 show that for a mean zero normal random variable $G$, for any $p > 0$, the random variable $|G|^p$ is $\alpha$-regular variation determining for any $\alpha > 0$. The $\alpha$-regular variation determinating property of gamma-type distributions was proved by an alternative Tauberian argument in [2] and [18].

EXAMPLE 4.5 (Pareto random variables and their reciprocals). For $p > 0$, a Pareto random variable $Y$ has density $f(x) = px^{-(p+1)}$ for $x > 1$. Note that $\log Y$ is exponentially distributed, hence infinitely divisible. By part (ii) of Corollary 4.3, we conclude that a Pareto random variable is $\alpha$-regular variation determining for any $\alpha < p$. This was also discovered in [20]. On the other hand, $1/Y$ (which has a particular beta distribution) is $\alpha$-regular variation determining for any $\alpha > 0$. This includes the standard uniform random variable.

EXAMPLE 4.6 (Lognormal random variables). By definition, for a lognormal random variable $Y$, $\log Y$ is normally distributed, hence infinitely divisible. By part (ii) of Corollary 4.3, we conclude that a lognormal random variable is $\alpha$-regular variation determining for any $\alpha > 0$.

EXAMPLE 4.7 (Cauchy random variables). Let $Y$ be the absolute value of the standard Cauchy random variable. Its logarithm has density $g(y) = ce^y(1 + e^{2y})^{-1}$ for $y > 0$ ($c$ is a normalizing constant). Therefore, $\log Y$ is an exponential tilting of a logistic random variable. Since the latter is infinitely divisible [30], so is $\log Y$. By part (ii) of Corollary 4.3, we conclude that the absolute value of a Cauchy random variable is $\alpha$-regular variation determining for $\alpha < 1$.

EXAMPLE 4.8 (A non-$\alpha$-regular variation determining random variable with truncated power tails). Let $0 < a < b < \infty$. For $\alpha > 0$, consider a random variable $Y$ with density $f(x) = cx^{-(\alpha+1)}, a < x < b$ ($c$ is a normalizing constant). Note that

$$E[Y^{\alpha+i\theta}] = c(e^{i\theta e^b} - e^{i\theta e^a}).$$

This last expression vanishes for

$$\theta = (e^b - e^a)^{-1} 2\pi n$$

for any $n \in \mathbb{Z}$. By Theorem 4.2, $Y$ is not $\alpha$-regular variation determining.

If a random variable $Y$ has an atom, then Corollary 2.2 applies, which gives the following result.

COROLLARY 4.9. Let $\alpha > 0$ and $Y$ be a positive random variable satisfying $EY^{\alpha+\delta} < \infty$ for some $\delta > 0$. Suppose that for some $x_0 > 0$, $P(Y = x_0) > 0$.

(i) If

$$(4.5) \qquad E[Y^\alpha \mathbf{1}(Y \neq x_0)] < x_0^\alpha P(Y = x_0)$$

then $Y$ is $\alpha$-regular variation determining.

(ii) If

$$(4.6) \qquad E[Y^\alpha \mathbf{1}(Y \neq x_0)] = x_0^\alpha P(Y = x_0),$$

then $Y$ is not $\alpha$-regular variation determining if and only if for some $\theta_0 > 0$,

$$(4.7) \qquad P[\log Y \in \log x_0 + 2\pi \theta_0 (2\mathbb{Z} + 1) | Y \neq x_0] = 1.$$

Example 4.10 (A two-point distribution). Let $a, b$ be two *different* positive numbers, and $Y$ a random variable taking values in the set $\{a, b\}$. If $a^\alpha P(Y = a) \neq b^\alpha P(Y = b)$, then part (i) of Corollary 4.9 applies, and $Y$ is $\alpha$-regular variation determining. If $a^\alpha P(Y = a) = b^\alpha P(Y = b)$, then part (ii) of Corollary 4.9 applies. Condition (4.7) is satisfied with, for example, $\theta_0 = |\log b - \log a|/(2\pi)$, and so $Y$ is not $\alpha$-regular variation determining.

Finally, let us consider the case $\alpha = 0$. We have the following simple statement.

Proposition 4.11. *Let $Y$ be a positive random variable independent of a random variable $Z$, such that $X = YZ$ has a slowly varying tail. If $EY^\delta < \infty$ for some $\delta > 0$, then $Z$ has a slowly varying tail, and $P(Z > x) \sim P(X > x)$ as $x \to \infty$.*

Proof. For any $\varepsilon > 0$,

$$P(X > x) \geq P(Y > \varepsilon)P(Z > x/\varepsilon)$$

and using slow variation of the tail of $X$ and letting $\varepsilon \to 0$ we conclude that

$$\limsup_{x \to \infty} \frac{P(Z > x)}{P(X > x)} \leq 1.$$

In particular, for some $C > 0$, $P(Z > x) \leq CP(X > x)$ for all $x > 0$. Therefore, for any $M > 0$, we can use Potter's bounds (see Proposition 0.8(ii) of [25]) to conclude that for all $x$ large enough

$$\frac{1}{P(X > x)} \int_M^\infty P(Z > x/y) F_Y(dy)$$

$$\leq C \int_M^\infty \frac{P(X > x/y)}{P(X > x)} F_Y(dy)$$

$$\leq C(1 + \delta) \int_M^\infty y^\delta F_Y(dy) \to 0, \qquad M \to \infty.$$

Since we also have

$$\frac{1}{P(X > x)} \int_0^M P(Z > x/y) F_Y(dy)$$

$$\leq P(Y \leq M) \frac{P(Z > x/M)}{P(X > x)}$$

$$= P(Y \leq M) \frac{P(X > x/M)}{P(X > x)} \frac{P(Z > x/M)}{P(X > x/M)},$$

we can use slow variation of the tail of $X$, and letting $M \to \infty$ to obtain

$$\liminf_{x \to \infty} \frac{P(Z > x)}{P(X > x)} \geq 1.$$

This completes the proof.  □

**5. The inverse problem for stochastic integrals.**  In this section, we are concerned with stochastic integrals with respect to Lévy processes as discussed in Example 1.3. We study the extent to which regular variation of the tail of the integral implies the corresponding regular variation of the tail of the Lévy measure of the Lévy process. Once again, the case $\alpha = 0$ is simple, and will be considered at the end of the section. For now, we assume that $\alpha > 0$.

THEOREM 5.1.  *Let $\alpha > 0$ and $(M(s))_{s \in \mathbb{R}}$ be a Lévy process with Lévy measure $\eta$. Let $f : \mathbb{R} \to \mathbb{R}_+$ be a measurable function such that for some $0 < \delta < \alpha$*

$$(5.1) \qquad \begin{cases} \displaystyle\int_{\mathbb{R}} [f(x)]^{\alpha - \delta} \vee [f(x)]^2 \, dx < \infty, & \text{if } 0 < \alpha < 2, \\ \displaystyle\int_{\mathbb{R}} [f(x)]^{\alpha - \delta} \vee [f(x)]^{\alpha + \delta} \, dx < \infty, & \text{if } \alpha \geq 2. \end{cases}$$

(i) *Assume that the integral $X = \int_{\mathbb{R}} f(s) M(ds)$ is well defined, and that the tail of $X$ is regularly varying with exponent $\alpha$. If*

$$(5.2) \qquad \int_{\mathbb{R}} [f(x)]^{\alpha + i\theta} \, dx \neq 0 \qquad \text{for all } \theta \in \mathbb{R},$$

*then the tail of the Lévy measure $\eta$ is regularly varying and (1.6) holds.*

(ii) *Let $f : \mathbb{R} \to \mathbb{R}_+$ be a measurable function satisfying (5.1). If (5.2) fails to hold, then there exists a Lévy process $(M(s))_{s \in \mathbb{R}}$ with Lévy measure $\eta$ which does not have a regularly varying tail, but the integral $X = \int_{\mathbb{R}} f(s) M(ds)$ is well defined and the tail of $X$ is regularly varying with exponent $\alpha$.*

PROOF.  If the integral defining the random variable $X$ is well defined, then $X$ is an infinitely divisible random variable with Lévy measure

$$(5.3) \qquad \eta_X = (\text{Leb} \times \eta) \circ T_f^{-1},$$

where $T_f : \mathbb{R} \times \mathbb{R} \to \mathbb{R}$ is given by $T_f(s, x) = x f(s)$; see [23]. Note that (5.3) says that $\eta_X = \eta \circledast \rho$, where $\rho = \text{Leb} \circ f^{-1}$ is a $\sigma$-finite measure on $\mathbb{R}_+$.

(i) The regular variation of the tail of $X$ is equivalent to the regular variation of the tail of $\eta_X$, in which case one also has

$$(5.4) \qquad \lim_{x \to \infty} \frac{P(X > x)}{\eta_X(x, \infty)} = 1,$$

see [14]. If the tail of $X$ is regularly varying, then we conclude that $\eta \circledast \rho$ has a regularly varying tail with exponent $\alpha$. The claim of part (i) of the theorem will follow once we check the conditions of part (i) of Theorem 2.3. The assumption (2.3) follows from (5.1), so one only has to verify the condition (2.18). Since $\eta$ is a Lévy measure, $\int_0^a w^2 \eta(dw) < \infty$ for every $0 < a < \infty$. Letting $p = \max(\alpha + \delta, 2) > \alpha$, (5.1) implies that for every $b, x > 0$

$$\int_0^b \rho(x/w, \infty) \eta(dw) \leq \int_{\mathbb{R}} |f(s)|^p \, ds \, x^{-p} \int_0^b w^p \eta(dw)$$
$$= o((\eta \circledast \rho)(x, \infty)) \qquad \text{as } x \to \infty.$$

This verifies (2.18), and hence completes the proof of part (i) of the theorem.

(ii) We use the corresponding part of Theorem 2.3. The first step is to check that the measure $\nu$ constructed as a counterexample there can be chosen to be a Lévy measure. In fact, this measure can be chosen to be a finite measure, as the construction in the proof of Theorem 4.2 shows, with a virtually identical argument. Call this Lévy measure $\eta_1$; it is a measure on $(0, \infty)$. Let $\eta = \eta_1 + \eta_1(-\cdot)$; then $\eta$ is a symmetric Lévy measure on $\mathbb{R}$. By construction, $\eta$ does not have a regularly varying tail.

Let $M$ be a symmetric Lévy process without a Gaussian component with Lévy measure $\eta$. Next, we check that the integral $X = \int_{\mathbb{R}} f(s) M(ds)$ is well defined. This is equivalent to verifying the condition

$$(5.5) \qquad \int_{\mathbb{R}} \int_{\mathbb{R}} [xf(s)]^2 \wedge 1 \, \eta(dx) \, ds < \infty,$$

see [23]. Rewrite the integral as

$$(5.6) \qquad \begin{aligned} &\int_{\mathbb{R}} \text{Leb}(\{s : f(s) > 1/|x|\}) \eta(dx) \\ &\quad + \int_{\mathbb{R}} x^2 \int_{\mathbb{R}} [f(s)]^2 \mathbf{1}(f(s) \leq 1/|x|) \, dx \, \eta(dx). \end{aligned}$$

Let $p = \max(2, \alpha + \delta)$. We can bound the first term in the right-hand side of (5.6) by

$$\int_{\mathbb{R}} (|x|^{-(\alpha-\delta)} \wedge |x|^{-p}) \eta(dx) \int_{\mathbb{R}} [f(s)]^{\alpha-\delta} \vee [f(s)]^p \, ds < \infty$$

in view of (5.1) and because $\eta$ is a Lévy measure. Consider the second term in the right-hand side of (5.6). The integral over the set $\{|x| \leq 1\}$ is finite since $g \in L^2$ and $\eta$ is a Lévy measure. Finally, the integral over the set $\{|x| > 1\}$ is bounded by

$$\int_{\mathbb{R}} x^2 \mathbf{1}(|x| > 1) |x|^{-(2-\alpha+\delta)} \eta(dx) \int_{\mathbb{R}} [f(s)]^{\alpha-\delta} \, ds < \infty,$$

because $f \in L^{\alpha - \delta}$, and the tail of the Lévy measure $\eta$ satisfies $\eta(x, \infty) \leq cx^{-\alpha}$ for all $x > 0$ and some finite constant $c$. Therefore, (5.5) holds. Hence, the integral $X = \int_{\mathbb{R}} f(s)M(ds)$ is well defined. It follows that the Lévy measure of the integral $X$ is given by (5.3), which by its construction has a regularly varying tail. This completes the argument. □

EXAMPLE 5.2. The *Ornstein–Uhlenbeck process with respect to a Lévy motion* (or a Lévy process) is defined by using the kernel $f(s) = e^{-\lambda s} \mathbf{1}(s > 0)$. With this kernel the left-hand side of (5.2) is equal to $\lambda(1 + i\theta)$ which does not vanish for real $\theta$. Therefore, if the marginal distribution of an Ornstein–Uhlenbeck process is regularly varying, the same holds for the Lévy measure of the underlying Lévy process. The same is true for the double-sided Ornstein–Uhlenbeck process, for which $f(s) = e^{-\lambda|s|}$. The Ornstein–Uhlenbeck process with respect to a symmetric $\alpha$-stable Lévy motion was already considered in Section 5 of [12]. Ornstein–Uhlenbeck processes with respect to a Lévy motion are well studied; see, for example, [28].

If there is a set of positive measure on which the function $f$ is constant, then Corollary 2.2 applies. The result parallels Corollary 4.9.

COROLLARY 5.3. *Suppose that for some $a > 0$, $\mathrm{Leb}(\{s : f(s) = a\}) > 0$.*

(i) *Suppose that*

$$\int_{\mathbb{R}} [f(s)]^{\alpha} \mathbf{1}(f(s) \neq a) \, ds < a^{\alpha} \, \mathrm{Leb}(\{s : f(s) = a\}).$$

*Then regular variation with exponent $\alpha$ of the tail of $X$ implies regular variation of the tail of the Lévy measure of the Lévy process.*

(ii) *If*

$$\int_{\mathbb{R}} [f(s)]^{\alpha} \mathbf{1}(f(s) \neq a) \, ds = a^{\alpha} \, \mathrm{Leb}(\{s : f(s) = a\}),$$

*then regular variation with exponent $\alpha$ of the tail of $X$ fails to imply regular variation of the tail of the Lévy measure of the Lévy process if and only if for some $\theta_0 > 0$, $\mathrm{Leb}(S_{a_0}^c) = 0$, where $S_{a_0}$ is defined in (2.16).*

A result which is more precise about the Lévy measure of the driving Lévy process than Theorem 5.1 can be obtained when the stochastic integral is infinite variance stable, see Theorem 6.1 below.

Finally, we consider the case $\alpha = 0$, where we have the following statement.

PROPOSITION 5.4. *Let $(M(s))_{s \in \mathbb{R}}$ be a Lévy process with Lévy measure $\eta$ and $f : \mathbb{R} \to \mathbb{R}_+$ be a measurable function such that $\mathrm{Leb}(\{s \in \mathbb{R} : f(s) \neq 0\}) < \infty$ and $\int_{\mathbb{R}} [f(s)]^2 \, ds < \infty$. Assume that the integral $X = \int_{\mathbb{R}} f(s)M(ds)$ is well defined, and that the tail of $X$ is slowly varying. Then the tail of the Lévy measure $\eta$ is slowly varying and (1.6) holds.*

PROOF. Since $X$ is well defined, and its tail is slowly varying relations (5.3) and (5.4) hold. As in the proof of Theorem 5.1, let $\rho$ be the image on $(0, \infty)$ of Lebesgue measure under the measurable map $g$. In the present case, $\rho$ is a finite measure. We have for $\varepsilon > 0$,

$$\eta_X(x, \infty) \geq \int_\varepsilon^\infty \eta(x/y, \infty) \rho(dy) \geq \eta(x/\varepsilon, \infty) \rho(\varepsilon, \infty).$$

Now using (5.4), slow variation and letting $\varepsilon \to 0$, we conclude that

(5.7) $$\limsup_{x \to \infty} \frac{\eta(x, \infty)}{P(X > x)} \leq \frac{1}{\text{Leb}(\{s \in \mathbb{R} : f(s) \neq 0\})}.$$

As in the proof of Proposition 4.11, the matching lower bound will follow once we show that

$$\lim_{M \to \infty} \limsup_{x \to \infty} \frac{1}{P(X > x)} \int_M^\infty \eta(x/y, \infty) \rho(dy) = 0.$$

To this end, for fixed $M > 0$ and $x > M$ write

$$\frac{1}{P(X > x)} \int_M^\infty \eta(x/y, \infty) \rho(dy) = \frac{1}{P(X > x)} \left( \int_M^x + \int_x^\infty \right) \eta(x/y, \infty) \rho(dy).$$

The first term in the right-hand side above can be bounded by (5.7) and Potter's bounds as

$$\int_M^x \frac{\eta(x/y, \infty)}{P(X > x/y)} \frac{P(X > x/y)}{P(X > x)} \rho(dy) \leq C \int_M^\infty y^2 \rho(dy)$$

for some $C > 0$. The right-hand side converges to zero as $M \to \infty$. Similarly, the second term above can be bounded by the fact that $\eta$ is a Lévy measure. Hence,

$$\frac{1}{P(X > x)} \int_x^\infty C x^{-2} y^2 \rho(dy),$$

which converges to zero as $x \to \infty$ because $X$ has a slowly varying tail. This completes the proof. □

## 6. Some identification problems for stable laws.

Non-Gaussian $\alpha$-stable distributions are known to have regularly varying tail with exponent $\alpha \in (0, 2)$, see [17] for the definition and properties of stable distributions and [27] for the case of general stable processes. In addition, they enjoy the property of convolution closure. Consider, for example, a sequence $(Z_i)$ of i.i.d. strictly $\alpha$-stable random variables for some $\alpha \in (0, 2)$. In this case,

$$X = \sum_{j=1}^\infty \psi_j Z_j \stackrel{d}{=} Z \left( \sum_{j=1}^\infty \psi_j^\alpha \right)^{1/\alpha}$$

provided the right-hand side is finite. Moreover, $P(|X| > x) \sim cx^{-\alpha}$ for some positive $c$ and the limits

$$\lim_{x \to \infty} P(X > x)/P(|X| > x) \quad \text{and} \quad \lim_{x \to \infty} P(X \leq -x)/P(|X| > x)$$

exist; see [17], Theorem XVII.5.1.

Similarly, if $(M(s))_{s \in \mathbb{R}}$ is symmetric $\alpha$-stable Lévy motion then

$$\int_{\mathbb{R}} f \, dM \overset{d}{=} M(1) \left( \int_{\mathbb{R}} |f(s)|^{\alpha} \, ds \right)^{1/\alpha},$$

where the left-hand side is well defined if and only if the integral on the right-hand side is finite for any real-valued function $f$ on $\mathbb{R}$; see [27], Chapter 3.

These examples show that the inverse problems in the case of stable distributions can be formulated in a more precise way: Given that the output of a linear filter is $\alpha$-stable is the input also $\alpha$-stable? We will not pursue here this question in greatest possible generality, but rather give an answer in two situations, in order to illustrate a different application of the cancellation property of $\sigma$-finite measures.

We start by considering stable stochastic integrals.

PROPOSITION 6.1. *Let $(M(s))_{s \in \mathbb{R}}$ be a symmetric Lévy process without Gaussian component. Suppose that for some measurable function $f : \mathbb{R} \to \mathbb{R}$ the stochastic integral $\int_{\mathbb{R}} f(s) \, dM(s)$ is well defined and represents a symmetric $\alpha$-stable random variable, where $0 < \alpha < 2$. If*

$$(6.1) \qquad \int_{\mathbb{R}} |f(s)|^{\alpha - \delta} \vee |f(s)|^{\alpha + \delta} \, ds < \infty \qquad \text{for some } 0 < \delta < \alpha$$

*and*

$$(6.2) \qquad \int_{\mathbb{R}} |f(s)|^{\alpha + i\theta} \, ds \neq 0 \qquad \text{for all } \theta \in \mathbb{R},$$

*then $(M(s))_{s \in \mathbb{R}}$ is an $\alpha$-stable Lévy process. Conversely, for any function $f$ which satisfies (6.1) but fails (6.2), there exists a nonstable symmetric Lévy process $(M(s))_{s \in \mathbb{R}}$ such that $\int_{\mathbb{R}} f(s) \, dM(s)$ has a symmetric $\alpha$-stable distribution.*

PROOF. Since $(M(s))_{s \in \mathbb{R}}$ is symmetric, the distributions of $\int_{\mathbb{R}} f \, dM$ and $\int_{\mathbb{R}} |f| \, dM$ are the same. Therefore, without loss of generality, we may assume that $f \geq 0$. Let $\nu$ be the Lévy measure of the law of $M(1)$ and $\rho$ denote the measure on $(0, \infty)$ induced by $f$, that is,

$$(6.3) \qquad \rho(B) = \int_{\mathbb{R}} \mathbf{1}_B(f(s)) \, ds, \qquad B \text{ a Borel set of } (0, \infty).$$

The stochastic integral $\int_{\mathbb{R}} f(s)\,dM(s)$ exists if and only if $\int_0^\infty \int_{\mathbb{R}} (xy)^2 \wedge 1\,\nu(dx)\rho(dy) < \infty$, and the Lévy measure $Q$ of the law of $\int_{\mathbb{R}} f\,dM$ is given by

$$Q(B) = \int_0^\infty \int_{\mathbb{R}} \mathbf{1}_{B\setminus\{0\}}(xy)\nu(dx)\rho(dy), \qquad B \text{ a Borel set of } \mathbb{R},$$

see [23]. By assumption, $Q$ is the Lévy measure of a symmetric $\alpha$-stable distribution. Therefore, for some $c > 0$, $\nu^{(+)} \circledast \rho = c\nu_\alpha$, where $\nu^{(+)}$ is the restriction of $\nu$ to $(0, \infty)$. This gives

$$(c^{-1}\|\rho\|_\alpha)(\nu^{(+)} \circledast \rho) = \nu_\alpha \circledast \rho,$$

where $\|\rho\|_\alpha = \int_0^\infty y^\alpha \rho(dy) < \infty$. By virtue of conditions (6.1) and (6.2), we may apply Theorem 2.1 and, therefore, we get by cancellation $(c^{-1}\|\rho\|_\alpha)\nu^{(+)} = \nu_\alpha$. Since $\nu$ is symmetric, it is the Lévy measure of an $\alpha$-stable distribution.

Conversely, let $\nu$ be a symmetric measure determined by (2.5) and (2.6), where $\rho$ is given by (6.3). Since the function $g$ in (2.6) is bounded, $\nu$ is a Lévy measure and, by (6.1), $\int_{\mathbb{R}} f\,dM$ exists. It follows that the Lévy measure $Q$ of $\int_{\mathbb{R}} f\,dM$ satisfies $Q(x, \infty) = \|\rho\|_\alpha x^{-\alpha}$ for $x > 0$, and since it is symmetric, $Q$ is the Lévy measure of an $\alpha$-stable distribution. □

It follows from Example 3.8 that it is possible to construct step functions $f \geq 0$ taking on only three distinct positive values such that (6.2) fails. However, it is impossible to achieve this with nonnegative functions taking on only two positive values.

Now assume that $X = \sum_{j=1}^\infty \psi_j Z_j$ is a symmetric $\alpha$-stable random variable for some $\alpha \in (0, 2)$, $Z$ is symmetric and $\sum_{j=1}^\infty |\psi_j|^{\alpha-\delta} < \infty$ for some $\delta \in (0, \alpha)$. Moreover, if we assume that $Z$ has an infinitely divisible distribution it follows from Theorem 6.1 that $Z$ is symmetric $\alpha$-stable provided $\sum_{j=1}^\infty |\psi_j|^{\alpha+i\theta}$ does not vanish for any $\theta \in \mathbb{R}$. Indeed, choose the integrand $f(s) = \psi_j$ for $s \in (j, j+1]$, $j \in \mathbb{Z}$, and the corresponding Lévy process $M$ on $\mathbb{R}$ such that $M(1) \overset{d}{=} Z$, that is, $M$ inherits the Lévy measure of $Z_1$. Then Theorem 6.1 implies that $M$ is symmetric $\alpha$-stable.

In this context, the following question arises: is it possible to drop the a priori condition of infinite divisibility on $Z$? The proposition below gives a partial answer. It considers finite sums, and assumes that the distribution of the sum is a totally skewed to the right stable distribution. This means that the Lévy measure of the stable distribution is concentrated on the positive half-line; see [27]. We find the answer in the cases $\alpha \in (0, 1) \cup (1, 2)$. What happens in the case $\alpha = 1$ is still unclear.

PROPOSITION 6.2. *Assume that $\alpha \in (0, 1) \cup (1, 2)$. Consider a set of positive $\alpha$-regular variation determining coefficients $\psi_1, \ldots, \psi_q$ for some $q \geq 2$,*

*that is,*

$$\text{(6.4)} \qquad \sum_{j=1}^{q} \psi_j^{\alpha+i\theta} \neq 0 \qquad \text{for all } \theta \in \mathbb{R}.$$

*Let $Z_1, \ldots, Z_q$ be i.i.d. random variables such that*

$$\text{(6.5)} \quad X_q = \sum_{j=1}^{q} \psi_j Z_j \text{ has a totally skewed to the right } \alpha\text{-stable distribution.}$$

*Then $Z$ has a totally skewed to the right $\alpha$-stable distribution. Conversely, if condition (6.4) does not hold, then there exists $Z$ that does not have an $\alpha$-stable law but $X_q$ is $\alpha$-stable.*

PROOF.  Assume that condition (6.4) holds.

*The case $0 < \alpha < 1$.* By subtracting a constant if necessary we may and do assume that $X_q$ has a strictly stable law. Write

$$\ell(s) = Ee^{-sZ} \quad \text{and} \quad G(x) = -\log \ell(x), \qquad x \geq 0.$$

Note that $G(0) = 0$ and $G$ is a nondecreasing continuous function. Therefore, there exists a $\sigma$-finite measure $\nu$ on $(0, \infty)$ such that $G(x) = \nu(0, x]$ for all $x > 0$. By assumption (6.5), there exists $c > 0$ such that

$$\text{(6.6)} \qquad \sum_{j=1}^{q} G(\psi_j x) = cx^\alpha \qquad \text{for all } x > 0,$$

which implies that

$$\text{(6.7)} \qquad \nu \circledast \rho = c\nu_{-\alpha},$$

where $\rho = \sum_{j=1}^{q} \delta_{1/\psi_j}$. By Theorem 2.1, this implies that $\nu = C\nu_{-\alpha}$ for some $C > 0$, as required.

*The case $1 < \alpha < 2$.* By subtracting expectations, we may and do assume that the random variables $Z_j$ and $X_q$ have zero mean. By Proposition 1.2.12 in [27] the Laplace transforms of $X_q$, and hence $Z$ are well defined for each nonnegative value of the argument.

As above, $\ell$ denotes the Laplace transform of $Z$, but now define $G(x) = \log \ell(x)$, $x \geq 0$. Note that $G(0) = 0$. Further, $\ell$ is a continuous convex function on $[0, \infty)$, and by the assumption of a zero mean of $X_q$, its one-sided derivative at zero is equal to zero. Hence, $\ell$ is nondecreasing on $[0, \infty)$. Therefore, $G$ is also a continuous nondecreasing function on $[0, \infty)$ and, therefore, there is still a $\sigma$-finite measure $\nu$ on $(0, \infty)$ such that $G(x) = \nu(0, x)$ for all $x > 0$.

By assumption (6.5) and Proposition 1.2.12 in [27], the relation (6.6) must still hold for some $c > 0$, which once again implies (6.7). As before, we conclude by Theorem 2.1 that for some $C > 0$, $G(x) = Cx^\alpha$ for all $x \geq 0$, which tell us that

$$\ell(x) = e^{Cx^\alpha}, \qquad x \geq 0.$$

Since the fact that Laplace transforms of two (not necessarily nonnegative) random variables coincide on an interval of positive length implies that the random variables are equal in law (page 390 or Problem 30.1 in [3]), we conclude that $Z$ must have a zero mean totally right skewed $\alpha$-stable law, as required.

If condition (6.4) fails to hold, then a counterexample can be constructed as in Theorem 6.1 and the comment following it. □

**Acknowledgments.** The authors would like to thank the two anonymous referees for constructive remarks and for a reference to the paper [29].

M. JACOBSEN
DEPARTMENT OF MATHEMATICS
UNIVERSITY OF COPENHAGEN
UNIVERSITETSPARKEN 5
DK-2100 COPENHAGEN
DENMARK
E-MAIL: martin@math.ku.dk

T. MIKOSCH
LABORATORY OF ACTUARIAL MATHEMATICS
UNIVERSITY OF COPENHAGEN
UNIVERSITETSPARKEN 5
DK-2100 COPENHAGEN
DENMARK
E-MAIL: mikosch@math.ku.dk

J. Rosiński
Department of Mathematics
121 Ayres Hall
University of Tennessee
Knoxville, Tennessee 37996-1300
USA
E-mail: rosinski@math.utk.edu

G. Samorodnitsky
School of Operations Research
  and Information Engineering
Cornell University
220 Rhodes Hall
Ithaca, New York 14853
USA
E-mail: gennady@orie.cornell.edu